# Embedding of partially ordered topological spaces in Fell topological hyperspaces

Jinlu Li
Department of Mathematics
Shawnee State University
Portsmouth Ohio 45662. USA
jli@shawnee.edu

**Abstract.** In this paper, we prove a theorem about embedding of some partially ordered topological spaces in topological hyperspaces equipped with Fell topology. Then we give some examples to show that the map defining the embedding may not be continuous with respect to Vietoris topology or Hausdorff topology equipped on the hyperspaces.

## 1. Introduction and preliminaries

In this section, we recall the concepts and some properties of partially ordered topological spaces and topological hyperspaces. We adopt the related definitions and notions used in [1−2] and [4].

### 1.1. Partially ordered topological spaces

Let $(X, \tau)$ be a topological space equipped with a partial order $\leqslant$. For every $x \in X$, we denote the $\leqslant$-principal ideal and principal filter of $x$ in $(X, \leqslant)$, respectively, as follows

$$x^{\downarrow} = \{u \in X: u \leqslant x\} \quad \text{and} \quad x^{\uparrow} = \{u \in X: u \geqslant x\}.$$

Let $A$ be a nonempty subset of $X$, we write

$$A^{\downarrow} = \{u \in X: u \leqslant x, \text{ for some } x \in A\},$$

and

$$A^{\uparrow} = \{u \in X: u \geqslant x, \text{ for some } x \in A\}.$$

For any $A \subseteq X$, we have $A \subseteq A^{\downarrow} \cap A^{\uparrow}$. If $A^{\downarrow} = A$, then $A$ is said to be $\leqslant$-decreasing and if $A^{\uparrow} = A$, then $A$ is said to be $\leqslant$-increasing.

In this paper, we follow the definition in [2] saying that $(X, \tau, \leqslant)$ is a partially ordered topological space, whenever the graph of $\leqslant$ is a closed subset of $X \times X$ with the product topology.

The Proposition 2.1 in the book "*topology and order*" by Nachbin provides a necessary and sufficient condition for partially ordered topological spaces.

**Proposition 2.1** in [4]. *Let $(X, \tau)$ be a topological space equipped with a partial order $\leqslant$. Then the following two properties are equivalent*:

(i) $(X, \tau, \leqslant)$ *is a partially ordered topological space (i.e. $\leqslant$ is a closed subset in $X \times X$)*;
(ii) *For $x, y \in X$, whenever $x \leqslant y$ false, there exist disjoint neighborhoods $V$ of $x$ and $U$ of $y$ such that $V$ is $\leqslant$-increasing and $U$ is $\leqslant$-decreasing*;

*Moreover,* (ii) *implies the following* (iii) *and* (iii) *does not imply* (ii).

(iii) *For every $x \in X$, both $x^{\downarrow}$ and $x^{\uparrow}$ are $\tau$-closed.*

### 1.2. Topologies on the hyperspaces of partially ordered topological spaces

Let $(X, \tau, \leqslant)$ be a partially ordered topological space. Let $C(X)$ denote the collection of all $\tau$-closed subsets of $X$ and $C_0(X) = C(X) \setminus \{\emptyset\}$. Denote the set of all $\leqslant$-principal ideals in $X$ by $C^{\downarrow}(X) = \{x^{\downarrow}: x \in X\}$. By Part (iii) in Proposition 2.1 in [4], we have $C^{\downarrow}(X) \subseteq C(X)$. In this subsection, we recall three ways to topology $C(X)$, which are used in this paper.

**1.2.1. The Fell topology $\tau_F$ on $C(X)$.** For every $\tau$-open subset $O$ of $X$ and every $\tau$-compact subsets $D$ of $X$, we denote the following subsets of $C(X)$:

$$O^- = \{A \in C(X): A \cap O \neq \emptyset\} \quad \text{and} \quad (X \setminus D)^+ = \{A \in C(X): A \cap D = \emptyset\}.$$

The set of subsets $O^-$ for $O$ running through the collection of all $\tau$-open subsets and $(X \setminus D)^+$ for $D$ running through all $\tau$-compact subsets of $X$ form the base of a topology on $C(X)$, which is called the Fell topology on $C(X)$ and it is denoted by $\tau_F$. Therefore, in this paper, $(C(X), \tau_F)$ is called the Fell topological hyperspace induced by this topological space $(X, \tau)$. Then $(C_0(X), \tau_F)$ and $(C^{\downarrow}(X), \tau_F)$ are also topological spaces with relative topology $\tau_F$.

Let $(X, \tau)$ be a Hausdorff topological space. The Fell topological hyperspace $(C(X), \tau_F)$ has many useful properties (see Beer [1] and Beer and Ok [2]). We list two of them below

(F_1) $(C(X), \tau_F)$ is compact;
(F_2) $x \to \{x\}$ embeds $X$ in $C(X)$.

**1.2.2. The Vietoris topology $\tau_V$ on $C(X)$.** The Vietoris topology $\tau_V$ on $C(G)$ is induced by the following base

$$\{O^-: O \text{ is an } \tau\text{-open subset of } X\} \cup \{(X \setminus E)^+: E \text{ is a } \tau\text{-closed subset of } X\}.$$

Then $(C(X), \tau_V)$ is called the Vietoris topological hyperspace of $(X, \tau)$.

**1.2.3. The Hausdorff topology $\tau_H$ on $C(X)$.** In particular, $(X, \tau)$ is a metric space and the topology $\tau$ is induced by a metric $d$ on $X$. The Hausdorff metric $H$ on $C(X)$ is defined, for any distinct $A, B \in C(X)$ as

$$H(A, B) = \max\left\{\sup_{a \in A}\left(\inf_{b \in B} d(a, b)\right), \sup_{b \in B}\left(\inf_{a \in A} d(b, a)\right)\right\}.$$

Then the Hausdorff metric $H$ induces the Hausdorff topology on $C(X)$, denoted by $\tau_H$ and $(C(X), \tau_H)$ is called the Hausdorff topological hyperspace of $(X, \tau)$.

**1.3. The canonical map $x \to x^{\downarrow}$**

When $(X, \tau, \preccurlyeq)$ is a partially ordered Hausdorrf topological space, for every $x \in X$, $\{x\} \in C(X)$. The property (F$_2$) in subsection 1.2.1 states that the map $x \to \{x\}$ embeds $(X, \tau)$ in $(C(X), \tau_F)$. In this subsection, we recall the canonical map $x \to x^{\downarrow}$ and its properties.

In a partially ordered topological space $(X, \tau, \preccurlyeq)$, $C^{\downarrow}(X) \subseteq C(X)\setminus\{\emptyset\}$ and $(C^{\downarrow}(X), \tau_F)$ is a topological subspace of $(C(X), \tau_F)$ with respect to the relative Fell topology. Let $(C(X), \tau_F)$ be equipped with the natural partial order $\subseteq$ (inclusion ordering). Since $\preccurlyeq$ and $\subseteq$ are partial orders on $X$ and $C(X)$, respectively, the canonical map is an order-embedding map, which satisfies that, for $x, y \in X$,

E$_1$. $x \preccurlyeq y$ if and only if $x^{\downarrow} \subseteq y^{\downarrow}$;
E$_2$. $x = y$ if and only if $x^{\downarrow} = y^{\downarrow}$.

We follow the concepts from [2]. We say that the canonical map $x \to x^{\downarrow}$ topologically order-embeds $(X, \tau, \preccurlyeq)$ in $(C(X), \tau_F, \subseteq)$, whenever, in addition to E$_1$ and E$_2$, this map $x \to x^{\downarrow}$ is a homeomorphism between $(X, \tau)$ and $(C^{\downarrow}(X), \tau_F)$ (with relative Fell topology $\tau_F$). More precisely speaking, the following two conditions are satisfied

E$_3$. The canonical map $x \to x^{\downarrow}$ is continuous from $(X, \tau)$ to $(C^{\downarrow}(X), \tau_F)$;
E$_4$. The map $x^{\downarrow} \to x$ is continuous from $(C^{\downarrow}(X), \tau_F)$ to $(X, \tau)$.

The topologically order-embedding from partially ordered topological spaces to hyperspaces have been studied by many authors (see [2–3]). In [2], Beer and Ok proved an embedding theorem on Hausdorff topological ∧-semilattices.

**Theorem 13 in [2]**. *Let $(X, \tau, \preccurlyeq)$ be a locally compact and order-connected Hausdorff topological ∧-semilattice. Then the canonical map $x \to x^{\downarrow}$ topologically order-embeds $(X, \tau, \preccurlyeq)$ in $(C(X), \tau_F, \subseteq)$.*

Recall that in the theories of ordered optimization and ordered variational inequalities, the considered partial orders equipped on the underlying spaces may be neither ∧-semilattice, nor ∨- semilattice. With the motivation of further possible applications of such embedding theorems to ordered optimizations and ordered variational inequalities, in section 2, we prove an embedding theorem on some partially ordered topological spaces, in which the conditions of order-connection and semilattice ordering relation are not required.

As applications of Proposition 2.1 in [4], the following results are proved in [2].

**Proposition 4** in [2]. *Let $(X, \tau, \preccurlyeq)$ be a locally compact partially ordered Hausdorff topological space. Then the following conditions are* equivalent:

(i)     *X is locally compact*;
(ii)    $\subseteq$ *is closed in $C(X) \times C(X)$ with the Fell topology*;
(iii)   $\subseteq$ *is closed in $C_0(X) \times C_0(X)$ with the relative Fell topology*.

## 1.4. The discontinuity of the canonical map $x \to x^{\downarrow}$ with respect to Vietoris and Hausdorff topologies

In contrast to Fell topology on the hyperspaces, in [2], two examples of totally discrete and locally compact metric (topological) $\wedge$-semilattices are provided, in which the map $x \to x^{\downarrow}$ is not continuous from $(X, \tau, \preccurlyeq)$ to $(C^{\downarrow}(X), \tau_V, \subseteq)$ and not continuous from $(X, \tau, \preccurlyeq)$ to $(C^{\downarrow}(X), \tau_H, \subseteq)$, respectively.

In the last section of this paper, we give two more examples of partially ordered topological spaces, which satisfy all conditions in Theorem 13 in [2]: # 1. A locally compact and order-connected metric (topological) lattice, in which the map $x \to x^{\downarrow}$ is not continuous at every point from $(X, \tau, \preccurlyeq)$ to $(C^{\downarrow}(X), \tau_V, \subseteq)$; #2. A locally compact and order-connected metric (topological) $\wedge$-semilattice, in which the map $x \to x^{\downarrow}$ is not continuous at every point (excepting one point) from $(X, \tau)$ to $(C^{\downarrow}(X), \tau_V)$ and not continuous at every point from $(X, \tau)$ to $(C^{\downarrow}(X), \tau_H)$.

## 2. Embedding partially ordered topological spaces in Fell topological hyperspaces

In this section, we introduce some concepts with respect to partially ordered topological spaces, which are used to prove the topologically embedding property of the order-embedding canonical map.

**Definition 2.1**. Let $(X, \tau, \preccurlyeq)$ be a partially ordered topological space. If, for any $x \in X$ and for any open subset $O$ in $X$, $x^{\downarrow} \cap O \neq \emptyset$ implies that there is an open neighborhood $D$ of $x$ such that $y^{\downarrow} \cap O \neq \emptyset$, for any $y \in D$, then $\preccurlyeq$ is said to be decreasing continuous on $(X, \tau)$.

**Proposition 2.2.** *Let $(X, \tau, \preccurlyeq)$ be a partially ordered Hausdorff topological space. If $\preccurlyeq$ is decreasing continuous on $(X, \tau)$, then the canonical map $x \to x^{\downarrow}$ from $(X, \tau)$ to $(C^{\downarrow}(X), \tau_F)$ is continuous.*

*Proof*. Since $x^{\downarrow}$ is $\tau$-closed, for any $x \in X$, it follows that the map $x \to x^{\downarrow}$ from $(X, \tau)$ into

$(C(X), \tau_F)$ is well-defined. With respect to the basis of $\tau_F$, we only prove the following two cases.

Case 1. Take an arbitrary point $x \in X$ with an arbitrary $\tau_F$-neighborhood $O^-$ of $x^\downarrow$, where $O$ is an nonempty open subset of $X$. It implies that $x^\downarrow \cap O \neq \emptyset$. Since $\leqslant$ is decreasing continuous on $(X, \tau)$, from $x^\downarrow \cap O \neq \emptyset$, there is an open neighborhood $D$ of $x$ such that $y^\downarrow \cap O \neq \emptyset$, for any $y \in D$. That is, $y^\downarrow \in O^-$, for any $y \in D$.

Case 2. Take a nonempty compact subset $E$ of $X$ with $x^\downarrow \cap E = \emptyset$, that is, $x^\downarrow \in (X \backslash E)^+$, and therefore, $(X \backslash E)^+$ is an $\tau_F$-neighborhood of $x^\downarrow$. The proof of this case is exactly same with the second part of the proof of Theorem 4.3 in [2]; and therefore, it is omitted here. □

**Definition 2.3.** Let $(X, \tau, \leqslant)$ be a partially ordered topological space.

(I) (Proper inclusion property) If, for any $x, y \in X$ with $y \leqslant x$,

$$y \in \text{int } x^\downarrow \Longrightarrow y^\downarrow \subseteq \text{int } x^\downarrow \quad \text{and} \quad x \in \text{int } y^\uparrow \Longrightarrow x^\uparrow \subseteq \text{int } y^\uparrow,$$

then $\leqslant$ is said to have proper inclusion property;

(II) (Dense boundaries) If, for any $x \in X$ and any open neighborhood $O$ of $x$, there are $a$, $b \in X$ with $b \leqslant a$ and $x \in \text{int } (b^\uparrow \cap a^\downarrow)$ such that, for any $y \in X$ with $y \notin b^\uparrow \cap a^\downarrow$ and for any $z \in \text{int } (b^\uparrow \cap a^\downarrow)$, we have that

(a) $y \leqslant z$ implies that there is $u \in \partial(b^\uparrow \cap a^\downarrow)$ such that $y \leqslant u \leqslant z$;
(b) $z \leqslant y$ implies that there is $v \in \partial(b^\uparrow \cap a^\downarrow)$ such that $z \leqslant v \leqslant y$,

then, $(X, \tau, \leqslant)$ is said to have dense boundaries.

(III) (Singular point) Let $x \in X$. If, for any open neighborhood $O$ of $x$, there is an open neighborhood $O_1$ of $x$ with $O_1 \subseteq O$, such that $\{y \in X \backslash O_1 : u \leqslant y, \text{ for some } u \in O_1\} = \emptyset$, then, $x$ is called an $\leqslant$-upper singular point.

(IV) (Upper compact bounded point) Let $x \in X$. If, for any open neighborhood $O$ of $x$, there are points $a, b \in O$ with $b \leqslant a$ such that

$$b^\uparrow \cap a^\downarrow \subseteq O, \; x \in \text{int}(b^\uparrow \cap a^\downarrow)$$

and $(b^\uparrow \cap a^\downarrow) \backslash (\text{int } a^\downarrow)$ is nonempty and compact,

then, $x$ is called an $\leqslant$-upper compact bounded point.

**Proposition 2.4.** *Let $(X, \tau, \leqslant)$ be a partially ordered Hausdorff topological space. Suppose that $(X, \tau, \leqslant)$ satisfies the following conditions*:

(i) *It has proper inclusion property*;

(ii) *It has dense boundaries;*
(iii) *For any point $x \in X$, $x$ is either an $\leqslant$-upper singular point, or an $\leqslant$-upper compact bounded point.*

*Then the map $x^{\downarrow} \to x$ is continuous from $(C^{\downarrow}(X), \tau_F)$ (with relative Fell topology) to $(X, \tau)$.*

*Proof.* Define a map $F$ from $C^{\downarrow}(X)$ to $X$ by

$$F(x^{\downarrow}) = x, \text{ for every } x^{\downarrow} \in C^{\downarrow}(X).$$

To prove that $F$ is continuous at every point $x^{\downarrow} \in C^{\downarrow}(X)$, based on the two cases in condition (iii) in this proposition, the proof is divided to two cases.

**Case 1.** For $x^{\downarrow} \in C(X)$, $x$ is an $\leqslant$-upper singular point. In this case, for any open neighborhood $O$ of $x$, there is an open neighborhood $O_1$ of $x$ with $O_1 \subseteq O$ such that

$$\{y \in X \setminus O_1 : u \leqslant y, \text{ for some } u \in O_1\} = \emptyset. \tag{2.1}$$

Since $x \in O_1$, it follows that $O_1^-$ is an $\tau_F$-open neighborhood of $x^{\downarrow}$ in $C(X)$. For all $y \in O_1$, it follows $y^{\downarrow} \in O_1^-$. We obtain

$$F(O_1^-) \supseteq O_1. \tag{2.2}$$

On the other hand, from (2.1), we have

$$y \notin O_1 \implies y^{\downarrow} \notin O_1^-, \text{ for } y \in X. \tag{2.3}$$

Combing (2.2) and (2.3), it follows that $F(O_1^-) = O_1 \subseteq O$.

**Case 2.** $x$ is an $\leqslant$-upper compact bounded point. In this case, for any open neighborhood $O$ of $x$, by conditions (iii), there are points $a, b \in O$ such that

$$b^{\uparrow} \cap a^{\downarrow} \subseteq O \subseteq E \text{ and } x \in \text{int}(b^{\uparrow} \cap a^{\downarrow}).$$

Let

$$O_1 = \text{int } (b^{\uparrow} \cap a^{\downarrow}) \subseteq O,$$

and

$$E_1 = (b^{\uparrow} \cap a^{\downarrow}) \setminus (\text{int } a^{\downarrow}).$$

Then $O_1$ is an open neighborhood of $x$, and $E_1$ is a nonempty and compact subset of $X$.

From $x \in O_1$, it follows that $x^{\downarrow} \in O_1^-$. On the other hand, by condition (i) and from $x \in \text{int } (b^{\uparrow} \cap a^{\downarrow}) \subseteq \text{int } a^{\downarrow}$, it implies that $x^{\downarrow} \subseteq \text{int } a^{\downarrow}$; and therefore, $x^{\downarrow} \cap E_1 = \emptyset$. That is $x^{\downarrow} \in (X \setminus E_1)^+$. We obtain

$$x^{\downarrow} \in O_1^- \cap (X \setminus E_1)^+.$$

Hence, $O_1^- \cap (X \setminus E_1)^+$ is an $\tau_F$-open neighborhood of $x^{\downarrow}$ in $C(X)$. In case 2, we show that

$$F(O_1^- \cap (X\backslash E_1)^+) = O_1. \tag{2.4}$$

At first, similar to the proof of $x^\downarrow \in O_1^- \cap (X\backslash E_1)^+$, we can show that,

$$y^\downarrow \in O_1^- \cap (X\backslash E_1)^+, \text{ for all } y \in O_1.$$

Hence, we obtain
$$F(O_1^- \cap (X\backslash E_1)^+) \supseteq O_1. \tag{2.5}$$

For the proof of (2.4), next we prove that, for any $y \in X$,

$$y \notin O_1 \implies y^\downarrow \notin O_1^- \cap (X\backslash E_1)^+. \tag{2.6}$$

For convenience, for any $S \subseteq X$, we denote $\bar{S} = X\backslash S$. Notice $O_1 \subseteq a^\downarrow$, we have

$$\overline{O_1} = (a^\downarrow \cup \overline{a^\downarrow}) \cap \overline{O_1} = (a^\downarrow \cap \overline{O_1}) \cup (\overline{a^\downarrow} \cap \overline{O_1}) = (a^\downarrow \cap \overline{O_1}) \cup \overline{a^\downarrow \cup O_1} = (a^\downarrow \backslash O_1) \cup \overline{a^\downarrow}.$$

So the proof of (2.6) is divided to the following two subcases.

**Subcase 2.1.** $y \in a^\downarrow \backslash O_1$. In this subcase, assume, by the way of contradiction, that $y^\downarrow \in O_1^-$, that is, $y^\downarrow \cap O_1 \neq \emptyset$. Then, there is $z \in y^\downarrow \cap O_1$ satisfying $b \leqslant z \leqslant y \leqslant a$.

Notice that $O_1 = \text{int}(b^\uparrow \cap a^\downarrow) = (\text{int } b^\uparrow) \cap (\text{int } a^\downarrow)$. From $z \in O_1$, it implies $z \in \text{int } b^\uparrow$. By condition (i), $z^\uparrow \subseteq \text{int } b^\uparrow$. By $z \leqslant y$, it follows that $y \in z^\uparrow \subseteq \text{int } b^\uparrow$. By the assumption $y \in a^\downarrow \backslash O_1$ in this case, it implies

$$\begin{aligned}
y &\in (\text{int } b^\uparrow) \cap (a^\downarrow \backslash O_1) \\
&= (\text{int } b^\uparrow) \cap (a^\downarrow \backslash (\text{int}(b^\uparrow \cap a^\downarrow))) \\
&= (\text{int } b^\uparrow) \cap (a^\downarrow \backslash ((\text{int } b^\uparrow) \cap (\text{int } a^\downarrow))) \\
&= (\text{int } b^\uparrow) \cap (a^\downarrow \cap \overline{(\text{int } b^\uparrow) \cap (\text{int } a^\downarrow)}) \\
&= (\text{int } b^\uparrow) \cap (a^\downarrow \cap (\overline{\text{int } b^\uparrow} \cup \overline{\text{int } a^\downarrow})) \\
&= ((\text{int } b^\uparrow) \cap a^\downarrow \cap \overline{\text{int } b^\uparrow}) \cup ((\text{int } b^\uparrow) \cap a^\downarrow \cap \overline{\text{int } a^\downarrow}) \\
&= \emptyset \cup ((\text{int } b^\uparrow) \cap a^\downarrow \cap \overline{\text{int } a^\downarrow}) \\
&= ((\text{int } b^\uparrow) \cap a^\downarrow \cap \overline{\text{int } a^\downarrow}) \\
&= ((\text{int } b^\uparrow) \cap a^\downarrow) \backslash (\text{int } a^\downarrow) \\
&\subseteq (b^\uparrow \cap a^\downarrow) \backslash (\text{int } a^\downarrow) \\
&= E_1.
\end{aligned}$$

We obtain $y \in E_1$. It follows that $y^\downarrow \notin (X\backslash E_1)^+$, which proves (2.6) in subcase 2.1.

**Subcase 2.2.** $y \in \overline{a^\downarrow}$, that is, $y \leqslant a$ false. In this case, if $y^\downarrow \cap O_1 = \emptyset$, then $y^\downarrow \notin O_1^-$, which proves (2.6). Suppose $y^\downarrow \cap O_1 \neq \emptyset$, that is, $y^\downarrow \in O_1^-$. Then, there is $z \in y^\downarrow \cap O_1$. By the assumption that

$y \leqslant a$ false, it implies $y \notin b^\uparrow \cap a^\downarrow$, and by condition (ii), there is $u \in \partial(b^\uparrow \cap a^\downarrow)$ such that

$$b \leqslant z \leqslant u \leqslant y.$$

It implies that
$$\begin{aligned}
u \in \partial(b^\uparrow \cap a^\downarrow) &= (b^\uparrow \cap a^\downarrow) \setminus \left(\text{int }(b^\uparrow \cap a^\downarrow)\right) \\
&= b^\uparrow \cap a^\downarrow \cap \overline{(\text{int }(b^\uparrow \cap a^\downarrow))} \\
&= (b^\uparrow \cap a^\downarrow) \cap \overline{((\text{int} b^\uparrow) \cap (\text{int} a^\downarrow))} \\
&= (b^\uparrow \cap a^\downarrow) \cap \overline{(\text{int } b^\uparrow \cup \text{ int } a^\downarrow)} \\
&= (b^\uparrow \cap a^\downarrow \cap \overline{\text{int } b^\uparrow}) \cup (b^\uparrow \cap a^\downarrow \cap \overline{\text{int } a^\downarrow}) \\
&= (b^\uparrow \cap a^\downarrow \cap \overline{\text{int } b^\uparrow}) \cup E_1.
\end{aligned} \qquad (2.7)$$

On the other hand, from $z \in y^\downarrow \cap O_1$ and $O_1 = \text{int }(b^\uparrow \cap a^\downarrow) \subseteq \text{int } b^\uparrow$, we have $z \in \text{int } b^\uparrow$. From condition (i), it follows that $z^\uparrow \subseteq \text{int } b^\uparrow$. By $b \leqslant z \leqslant u \leqslant a$, we have

$$u \in z^\uparrow \subseteq \text{int } b^\uparrow$$

It implies $u \notin b^\uparrow \cap a^\downarrow \cap \overline{\text{int } b^\uparrow}$. From (2.7), it follows that $u \in E_1$. From $b \leqslant z \leqslant u \leqslant y$, it implies that $y^\downarrow \cap E_1 \neq \emptyset$, that is $y^\downarrow \notin (X \setminus E_1)^+$, which proves (2.6) in subcase 2.2.

Hence (2.6) is proved, that is,
$$F(O_1^- \cap (X \setminus E_1)^+) \subseteq O_1. \qquad (2.8)$$

Then, (2.4) follows immediately from (2.5) and (2.8), which proves the continuity of $F$ at the arbitrary given point $x^\downarrow \in C^\downarrow(X)$. □

We provide a counter example below to show that condition (iii) in Proposition 2.4 is a necessary condition to assure the continuity of the map $x^\downarrow \to x$, from $(C^\downarrow(X), \tau_F)$ to $(X, \tau)$.

**Example 2.5.** Define
$$X = \{(u, v) \in \mathbb{R}^2 : uv > 0\} \cup \{\theta\}, \text{ where } \theta = (0, 0).$$

Let $X$ be equipped with the standard Euclidean topology $\tau$ and equipped with the coordinate wise partial order $\leqslant$. Then $(X, \tau, \leqslant)$ is a partially ordered topological space. Moreover,

(a) $\theta$ is neither an $\leqslant$-upper singular point, nor an $\leqslant$-upper compact bounded point;
(b) The map $x^\downarrow \to x$ from $(C^\downarrow(X), \tau_F)$ to $(X, \tau)$ is not continuous at $\theta^\downarrow$.

*Proof.* It is clear to see that $\theta$ is not an $\leqslant$-upper singular point. Then, we show that $\theta$ is not an $\leqslant$-upper compact bounded point. As a matter of fact, for any $\tau$-open neighborhood $O$ of $\theta$, and for any given two points $a, b \in O$ with $b \leqslant a$ satisfying $b^\uparrow \cap a^\downarrow \subseteq O$ and $x \in \text{int}(b^\uparrow \cap a^\downarrow)$ (One can see that the existence of such points). It follows that $b < \theta < a$. However, $(b^\uparrow \cap a^\downarrow) \setminus (\text{int } a^\downarrow)$ is nonempty and NOT compact. It implies that condition (iii) in Proposition 2.4 does not hold. One can prove that the map $x^\downarrow \to x$ from $(C^\downarrow(X), \tau_F)$ to $(X, \tau)$ is not continuous at $\theta^\downarrow$.

Combining Propositions 2.2 and 2.4, and from Proposition 1.1, we have the following embedding theorem.

**Theorem 2.6.** *Let $(X, \tau, \preccurlyeq)$ be a partially ordered Hausdorff topological space with $\preccurlyeq$ being decreasing continuous. Suppose that $(X, \tau, \preccurlyeq)$ satisfies the following conditions:*

(i)  *It has proper inclusion property;*
(ii) *It has dense boundaries;*
(iii) *For any point $x \in X$, $x$ is either an $\preccurlyeq$-upper singular point, or an $\preccurlyeq$-upper compact bounded point.*

*Then map $x \to x^{\downarrow}$ topologically order-embeds $(X, \tau, \preccurlyeq)$ in $(C(X), \tau_F, \subseteq)$.*

## 3. Applications to subsets of finite-dimensional partially ordered topological vector spaces

Based on our motivation of applying the embedding theory to variational analysis, in this section, we study the properties of hyperspaces of finite-dimensional partially ordered topological vector spaces. In [2], an embedding theorem in finite-dimensional spaces was proved.

**Corollary 29 in [2].** *Let $(X, \tau, \preccurlyeq)$ be a finite-dimensional partially ordered topological vector space. Then the canonical map $x \to x^{\downarrow}$ topologically order-embeds $(X, \tau, \preccurlyeq)$ in $(C(X), \tau_F, \subseteq)$.*

Let $(X, \tau, \preccurlyeq)$ be a given finite-dimensional partially ordered topological vector space. By the above Corollary 29 in [2], it satisfies that the canonical order-embedding is a topological embedding from $(X, \tau, \preccurlyeq)$ to $(C(X), \tau_F, \subseteq)$. If we take an arbitrary subset $Y \subseteq X$ satisfying that the subspace $(Y, \tau, \preccurlyeq)$ is a partially ordered topological space, in which $\tau$ is the relative topology restricted on $Y$, and the hyperspace space $(C(Y), \tau_F, \subseteq)$ is equipped with the relative Fell topology restricted on $C(Y)$. The canonical map $y \to y^{\downarrow}$ is always an order-embedding map from $(Y, \tau, \preccurlyeq)$ to $(C(Y), \tau_F, \subseteq)$. Then one may ask the following question:

(Q). Does the canonical map $y \to y^{\downarrow}$ topologically embed $(Y, \tau)$ in $(C(Y), \tau_F)$?

It is a little surprising to see that, in general, the answer is negative. Moreover, with the relative topology $\tau_F$ on $C^{\downarrow}(Y)$, we have

(a) The canonical map $y \to y^{\downarrow}$ from $(Y, \tau)$ to $(C^{\downarrow}(Y), \tau_F)$ may not be continuous;
(b) The map $y^{\downarrow} \to y$ from $(C^{\downarrow}(Y), \tau_F)$ to $(Y, \tau)$ may not be continuous.

Throughout this section, let $\mathbb{R}^n$ denote the *n*-d Euclidean space, for some $n > 1$, equipped with the standard topology $\tau$ induced by the ordinary norm $\|\cdot\|$. Let $K$ be a pointed closed and convex cone in $\mathbb{R}^n$ with nonempty interior, which induces a partial order $\preccurlyeq$ on $\mathbb{R}^n$. Then $(\mathbb{R}^n, \tau, \preccurlyeq)$ is a finite-dimensional partially ordered topological vector space. By Corollary 29 in [2], the canonical map $x \to x^{\downarrow}$ topologically order-embeds $(\mathbb{R}^n, \tau, \preccurlyeq)$ in $(C(\mathbb{R}^n), \tau_F, \subseteq)$. In this section, we apply Theorem 2.6 in last section and with some approaches different from the proof of Corollary 29 in [2] to exam the subsets of $\mathbb{R}^n$ for answering the question (Q) with respect to open and closed cases. We first prove the following three lemmas in $(\mathbb{R}^n, \tau, \leqslant)$.

**Lemma 3.1**. *For any $a \in K$, $a^{\downarrow} \cap (-a)^{\uparrow}$ is bounded.*

*Proof*. It is trivial for $a = \theta$. We suppose that $a \in K \setminus \{\theta\}$. Assume by contradiction that $a^{\downarrow} \cap (-a)^{\uparrow}$ is unbounded. There is a sequence $\{b_m\} \subseteq a^{\downarrow} \cap (-a)^{\uparrow}$ with $\|b_m\| = 1$, for $m = 1, 2, \ldots$, and there is an increasing sequence of positive numbers $\{t_m\}$ with limit infinity such that $\{t_m b_m\} \subseteq a^{\downarrow} \cap (-a)^{\uparrow}$. Suppose that $b_m \to b$, as $m \to \infty$. Then $\|b\| = 1$. Since $a^{\downarrow} \cap (-a)^{\uparrow}$ is closed, it follows that $b \in a^{\downarrow} \cap (-a)^{\uparrow}$. For any fixed $k > 1$, from $b_m \to b$, as $m \to \infty$, it implies that $t_k b_m \to t_k b$, as $m \to \infty$. The closeness of $a^{\downarrow} \cap (-a)^{\uparrow}$ implies that $t_k b \in a^{\downarrow} \cap (-a)^{\uparrow}$, for $k = 1, 2, \ldots$. Then

$$t_k b \leqslant a \text{ and } t_k b \geqslant -a, \text{ for } k = 1, 2, \ldots.$$

It follows that

$$b \leqslant \frac{a}{t_k} \text{ and } b \geqslant -\frac{a}{t_k}, \text{ for } k = 1, 2, \ldots.$$

Take $k \to \infty$, and by the closeness of $b^{\uparrow}$ and $b^{\downarrow}$, we have $b \leqslant \theta$ and $b \geqslant \theta$. It implies that $b = \theta$, which contradicts to $\|b\| = 1$. □

**Lemma 3.2**. *For any $x \in \mathbb{R}^n$ and any open neighborhood $U$ of $x$, there are points $a, b \in U$ with $b \leqslant x \leqslant a$, such that $x \in \text{int}(a^{\downarrow} \cap b^{\uparrow}) \subseteq a^{\downarrow} \cap b^{\uparrow} \subseteq U$.*

Proof. We only need to show the case that $x = \theta$ and $U = \{z \in \mathbb{R}^n : \|z\| < r\}$, for some $r > 0$. From the assumption that $\text{int}K \neq \emptyset$, it follows that $(\text{int}K) \cap U \neq \emptyset$. Take $a \in (\text{int}K) \cap U$. Then from $a \geqslant \theta$ and $\|a\| > 0$, we have

$$\theta \in \text{int}((ta)^{\downarrow} \cap (-ta)^{\uparrow}), \text{ for any } 0 < t < 1.$$

We claim that there is $0 < t_0 < 1$, such that

$$(t_0 a)^{\downarrow} \cap (-t_0 a)^{\uparrow} \subseteq U. \tag{3.1}$$

Assume by contradiction that (3.1) does not hold for any $0 < t < 1$. Then it follows that

$$(ta)^{\downarrow} \cap (-ta)^{\uparrow} \nsubseteq U, \text{ for any } 0 < t < 1.$$

That is,

$$((ta)^{\downarrow} \cap (-ta)^{\uparrow}) \cap (\mathbb{R}^n \setminus U) \neq \emptyset, \text{ for any } 0 < t < 1.$$

For any $0 < t < 1$, $((ta)^{\downarrow} \cap (-ta)^{\uparrow}) \cap (\mathbb{R}^n \setminus U)$ is closed, and by Lemma 3.1, it is bounded. Hence $((ta)^{\downarrow} \cap (-ta)^{\uparrow}) \cap (\mathbb{R}^n \setminus U)$ is compact. Since $s < t$ implies $(sa)^{\downarrow} \cap (-sa)^{\uparrow} \subseteq (ta)^{\downarrow} \cap (-ta)^{\uparrow}$, it follows that

$$\cap_{0 < t < 1} (((ta)^{\downarrow} \cap (-ta)^{\uparrow}) \cap (\mathbb{R}^n \setminus U)) \neq \emptyset.$$

Take $b \in \cap_{0 < t < 1}(((ta)^{\downarrow} \cap (-ta)^{\uparrow}) \cap (\mathbb{R}^n \setminus U))$. Then $\|b\| \geq r > 0$ satisfying

$$b \leqslant ta \text{ and } b \geqslant -ta, \text{ for any } 0 < t < 1.$$

Take $t \to 0$, and by the closeness of $b^\uparrow$ and $b^\downarrow$, we have $b \leqslant \theta$ and $b \geqslant \theta$. It implies that $b = \theta$, which contradicts to $\|b\| \geq r > 0$. □

**Lemma 3.3**. *For any $a, b, c \in \mathbb{R}^n$ with $b \leqslant a$, suppose that $c \notin a^\downarrow \cap b^\uparrow$. Then we have*

(i) *If $c \leqslant x$, for some $x \in \text{int}(a^\downarrow \cap b^\uparrow)$, then there is $y \in \partial(a^\downarrow \cap b^\uparrow)$ such that $c \leqslant y \leqslant x$;*

(ii) *If $x \leqslant c$, for some $x \in \text{int}(a^\downarrow \cap b^\uparrow)$, then there is $z \in \partial(a^\downarrow \cap b^\uparrow)$ such that $x \leqslant z \leqslant c$.*

*Proof of* (i). It is clear that $c \neq x$. So we suppose that $c < x$. Let $S(c, x)$ be the open segment in $\mathbb{R}^n$ with non-included ending points $c$ and $x$. From $c < x$, it implies that, for any $y \in S(c, x)$, we have $c \leqslant y \leqslant x$. Since $S(c, x) \cap (\text{int}(a^\downarrow \cap b^\uparrow)$ and $(c, x) \setminus (a^\downarrow \cap b^\uparrow)$ are two open subsets of $S(c, x)$, by the connectivity of $S(c, x)$, it implies that

$$S(c, x) \supsetneq (S(c, x) \cap (\text{int}(a^\downarrow \cap b^\uparrow))) \cup (S(c, x) \setminus (a^\downarrow \cap b^\uparrow))$$
$$= (S(c, x) \cap (\text{int}(a^\downarrow \cap b^\uparrow))) \cup (S(c, x) \cap \overline{a^\downarrow \cap b^\uparrow})$$
$$= S(c, x) \cap (\text{int}(a^\downarrow \cap b^\uparrow) \cup \overline{a^\downarrow \cap b^\uparrow})$$

It follows that $S(c, x) \cap (\overline{\text{int}(a^\downarrow \cap b^\uparrow)} \cap (a^\downarrow \cap b^\uparrow)) \neq \emptyset$. Noticing that $(\partial(a^\downarrow \cap b^\uparrow)) \cap S(c, x) = S(c, x) \cap (\overline{\text{int}(a^\downarrow \cap b^\uparrow)} \cap (a^\downarrow \cap b^\uparrow))$, we obtain

$$(\partial(a^\downarrow \cap b^\uparrow)) \cap S(c, x) \neq \emptyset.$$

Take a point $d \in (\partial(a^\downarrow \cap b^\uparrow)) \cap (c, x)$. Then $d \in \partial(a^\downarrow \cap b^\uparrow)$ and $c \leqslant d \leqslant x$. (ii) can be similarly proved. □

**Theorem 3.4**. *Let $Y$ be a nonempty open subset of $\mathbb{R}^n$ equipped with the relative topology and the partial order $\leqslant$ induced by $K$. Then $(Y, \tau, \leqslant)$ is a partially ordered topological space and the map $y \to y^\downarrow$ topologically order-embeds $(Y, \tau, \leqslant)$ in $(C(Y), \tau_F, \subseteq)$.*

*Proof.* By Proposition 2.1 in [4], the proof of the fact that $(Y, \tau, \leqslant)$ is a partially ordered topological space is straight forward and it is omitted here. Since $Y$ is open, it implies that $\leqslant$ is decreasing continuous on $(Y, \tau)$. From the assumption that $K$ has nonempty interior, one can see that $(Y, \tau, \leqslant)$ has $\leqslant$-proper inclusion property. Lemma 3.3 shows that $(Y, \tau, \leqslant)$ has dense boundaries. Finally, since $Y$ is an open subset in $\mathbb{R}^n$, by Lemma 3.2, it implies that every point $y \in Y$ is an $\leqslant$-upper compact bounded point. Hence, $(Y, \tau, \leqslant)$ satisfies all conditions in Theorem 2.6, which proves this theorem immediately. □

In contrast to the embedding property on the open subsets of $\mathbb{R}^n$, the situations of closed subsets are very different. We provide next example of closed and locally compact subset $Y$ of $\mathbb{R}^3$ such that the canonical map $y \to y^\downarrow$ from $(Y, \tau, \leqslant)$ to $(C^\downarrow(Y), \tau_F, \subseteq)$ is not continuous.

**Example 3.5**. Define two closed triangles $T_1$ and $T_2$ in $\mathbb{R}^3$ by

$T_1$ has vertexes $(0, 0, 0)$, $(0, -1, 0)$, $(-1, -1, 0)$;
$T_2$ has vertexes $(0, 0, 0)$, $(-1, -1, 0)$, $(-1, -1, -1)$.

Let $Y = T_1 \cup T_2$ be equipped with the standard Euclidean topology $\tau$ and the coordinate wise partial order $\leqslant$. Then

(i) $(Y, \tau, \leqslant)$ is a locally compact $\leqslant$-connected lattice (it is not a topological lattice);
(ii) $\wedge: Y \times Y \to Y$ is not continuous at point $((0, 0, 0), (0, -1, 0))$. It follows that $(Y, \tau, \leqslant)$ is not a topological $\wedge$-semilattice;
(iii) The canonical map $y \to y^{\downarrow}$ from $(Y, \tau, \leqslant)$ to $(C^{\downarrow}(Y), \tau_F, \subseteq)$ is not continuous at point $(0, 0, 0)$.

*Proof.* Proof of (i). It is clear to see that $(Y, \tau, \leqslant)$ is a closed and locally compact subset of $\mathbb{R}^3$. For any $y_1 = (u_1, v_1, w_1)$ and $y_2 = (u_2, v_2, w_2) \in Y$ with $y_1 \leqslant y_2$ and $y_1 \neq y_2$, we have

$$y_1^{\uparrow} \cap y_2^{\downarrow} = \{(u, v, w) \in Y: u_1 \leq u \leq u_2, v_1 \leq v \leq v_2 \text{ and } w_1 \leq w \leq w_2\}.$$

$y_1^{\uparrow} \cap y_2^{\downarrow}$ is a closed polygonal region in $\mathbb{R}^3$, which is indeed connected. It follows that $(Y, \tau, \leqslant)$ is $\leqslant$-connected. To check that $(Y, \leqslant)$ is a lattice, we first check that $(Y, \leqslant)$ is a $\vee$-semilattice. For arbitrary points $y_1 = (u_1, v_1, w_1)$ and $y_2 = (u_2, v_2, w_2) \in Y$, if $y_1, y_2 \in T_i$, for a fixed $i = 1, 2$. Then $(u_1 \vee u_2, v_1 \vee v_2, w_1 \vee w_2) \in T_i$ and such that $y_1 \vee y_2 = (u_1 \vee u_2, v_1 \vee v_2, w_1 \vee w_2)$. Next, suppose that $y_1$ and $y_2$ are in different $T_1$ and $T_2$. For example, $y_1 = (u_1, v_1, 0) \in T_1$ and $y_2 \in T_2$. One can check that $y_1 \vee y_2 = (u_1 \vee u_2, v_1 \vee v_2, 0)$. It follows that $(Y, \leqslant)$ is a $\vee$-semilattice. Moreover, $\vee Y = (0, 0, 0)$. We can similarly check that $(Y, \leqslant)$ is a $\wedge$-semilattice with $\wedge Y = (-1, -1, -1)$. Hence, $(Y, \leqslant)$ is a lattice.

Proof of (ii). Notice that $(0, 0, 0) \wedge (0, -1, 0) = (0, -1, 0)$. Take an arbitrary open neighborhood $G$ of $(0, -1, 0)$ with $G \subseteq T_1 \setminus T_2$. For any open neighborhoods $O$, $G_1$ of points $(0, 0, 0)$ and $(0, -1, 0)$, respectively, we have $O \cap (T_2 \setminus T_1) \neq \emptyset$. For any point $y = (u, v, w) \in O \cap (T_2 \setminus T_1)$, we have that $u, v, w < 0$. It follows that $y \wedge (0, -1, 0) = (u, -1, w) \in T_2 \setminus T_1$. By $G \subseteq T_1 \setminus T_2$, it follows that $y \wedge (0, -1, 0) \notin G$. It proves that $\wedge: Y \times Y \to Y$ is not continuous at point $((0, 0, 0), (0, -1, 0))$.

Proof of (iii). Take an arbitrary nonempty open subset $O \subseteq T_1$. One can see that $O \subseteq (0, 0, 0)^{\downarrow}$; and therefore, $O^+$ is an $\tau_F$-open neighborhood of $(0, 0, 0)^{\downarrow}$. On the other hand, for any $y = (u, v, w) \in T_2$ with $w < 0$, $y^{\downarrow} \cap O = \emptyset$, which implies $y^{\downarrow} \notin O^+$. It proves (iii).

Next example is a closed subset $Y$ of $\mathbb{R}^2$ that the map $y^{\downarrow} \to y$ from $(C^{\downarrow}(Y), \tau_F, \subseteq)$ to $(Y, \tau, \leqslant)$ is not continuous.

**Example 3.6.** Define
$$Y = \{(u, v) \in \mathbb{R}^2: u > 0, v \geq \frac{4}{u}\} \cup \{(u, v) \in \mathbb{R}^2: u^2 + v^2 \leq 1\}.$$

Let $Y$ be equipped with the standard Euclidean topology $\tau$ and equipped with the coordinate wise partial order $\preccurlyeq$. Then $(Y, \tau, \preccurlyeq)$ is a locally compact partially ordered topological space and the map $y^\downarrow \to y$ from $(C^\downarrow(Y), \tau_F, \subseteq)$ to $(Y, \tau, \preccurlyeq)$ is not continuous.

*Proof.* It is clear that there is no $\preccurlyeq$-upper singular point in $Y$. Let $U = \{(u, v) \in Y: u^2 + v^2 \leq 1\}$. Then, every point $(u, v)$ on the boundary of $U$ satisfying $0 \leq u, v \leq 1$, is not $\preccurlyeq$-upper compact bounded. Hence, $(Y, \tau, \preccurlyeq)$ does not satisfy condition (iii) in Theorem 2.6.

Define $F(y^\downarrow) = y$, for $y^\downarrow \in C^\downarrow(Y)$. We prove that $F$ is not continuous at point $(1,0)^\downarrow$. To this end, for arbitrary given $\tau$-open neighborhood $V \subseteq U$ of $(1,0)$, we show that, for some nonempty open subset $O$ and nonempty compact subset $E$, $O^+ \cap (Y \setminus E)^- \subseteq C^\downarrow(Y)$ is an $\tau_F$-open neighborhood of $(1,0)^\downarrow$, which satisfies the following property

$$F(O^+ \cap (Y\setminus E)^-) \nsubseteq V. \tag{3.2}$$

In (3.2), this $\tau_F$-open neighborhood $O^+ \cap (Y\setminus E)^-$ of $(1,0)^\downarrow$ is the "smallest type" elements in the basis of $\tau_F$ equipped on $C^\downarrow(Y)$. Based on our desire that $(1, 0) \in F(O^- \cap (Y\setminus E)^+)$, we must have $(1, 0) \in O$ and $(1, 0) \notin E$ satisfying $O \cap E = \emptyset$. Since $(1, 0) \in O$ and $O$ is open, there is $p > 0$ such that

$$O_1 \equiv \{(u, v) \in Y: u^2 + v^2 \leq 1, p < u \leq 1\} \subseteq O. \tag{3.3}$$

On the other hand, we have

$$(1,0)^\downarrow = \{(u, v) \in Y: u^2 + v^2 \leq 1, v \leq 0\}. \tag{3.4}$$

Since we desire that $(1,0)^\downarrow \in O_1^- \cap (Y\setminus E)^+$, that is, $(1,0)^\downarrow \cap E = \emptyset$. By $E \subseteq Y$ and (3.4), it follows that

$$E \subseteq \{(u, v) \in Y: v > 0\}.$$

Since $E$ is compact, it implies that there is $q > 0$, such that

$$E \cap \{(u, v) \in Y: 0 < v < q\} = \emptyset. \tag{3.5}$$

Take $d = \min\{\sqrt{1 - p^2}, q\}$. From (3.4) and (3.5), for any $0 < r < d$ and for any point $(u, v) \in Y$ with $0 < v < r$ and $u > \frac{4}{r}$, we have

$$(u, v)^\downarrow \cap O_1 \neq \emptyset, \text{ that is, } (u, v)^\downarrow \in O_1^-. \tag{3.6}$$

and

$$(u, v)^\downarrow \cap E = \emptyset, \text{ that is, } (u, v)^\downarrow \in (Y\backslash E)^+. \tag{3.7}$$

Combing (3.6) and (3.7), we have

$$(u, v)^\downarrow \in O_1^- \cap (Y\backslash E)^+, \text{ for any } (u, v) \in Y \text{ with } 0 < v < r < d \text{ and } u > \frac{4}{r}.$$

It follows that

$$F(O_1^- \cap (Y\backslash E)^+) \supset \{(u, v) \in Y: u > \frac{4}{r} \text{ and } 0 < v < r < d\}.$$

Since $\{(u, v) \in Y: u > \frac{4}{r} \text{ and } 0 < v < r < d\}$ is unbounded and $V$ is bounded, it proves (3.2). □

## 4. Vietoris and Hausdorff topologies on Hyperspaces of topological ∧-semilattices

Let $(X, \tau, \leqslant)$ be a partially ordered topological space. In addition, if $\leqslant$ is an ∧-semilattice and ∧: $X \times X \to X$ is a continuous operator, then $(X, \tau, \leqslant)$ is called a topological ∧-semilattice. We similarly define topological ∨-semilattices. If $(X, \tau, \leqslant)$ is both topological ∧-semilattice and topological ∨-semilattice, then $(X, \tau, \leqslant)$ is called a topological lattice.

Topological ∧-semilattices, topological ∨-semilattices and topological lattices all are considered as special cases of partially ordered topological spaces. In this section, we first recall an important embedding theorem on topological ∧-semilattices proved in [2].

**Theorem 13 in [2]**. *Let $(X, \tau, \leqslant)$ be a locally compact and order-connected Hausdorff topological ∧-semilattice. Then the canonical map $x \to x^\downarrow$ topologically order-embeds $(X, \tau, \leqslant)$ in $(C(X), \tau_F, \subseteq)$.*

In Example 3.5 in last section, $(Y, \tau, \leqslant)$ is a closed and locally compact subset of $\mathbb{R}^3$ and it is an $\leqslant$-connected lattice, in which the operator ∧: $Y \times Y \to Y$ is not continuous at point $((0, 0, 0), (0, -1, 0))$; and therefore, $(Y, \tau, \leqslant)$ is not a topological ∧-semilattice. Meanwhile, we showed that the canonical map $y \to y^\downarrow$ from $(Y, \tau, \leqslant)$ to $(C^\downarrow(Y), \tau_F, \subseteq)$ is not continuous at point $(0, 0, 0)$. Hence, Example 3.5 shows that the condition that the operator ∧: $Y \times Y \to Y$ is continuous in Theorem 13 in [2] is indeed a necessary condition to guarantee that the order-embedding $x \to x^\downarrow$ can topologically embeds $(Y, \tau, \leqslant)$ in $(C^\downarrow(Y), \tau_F, \subseteq)$.

In order to investigate the connections between Fell topology and vietoris topology or Hausdorff topology on hyperspaces, in [2], Beer and Ok studied the cases that if the hyperspaces of Hausdorff topological ∧-semilattices (as special cases of partially ordered topological spaces) are equipped with vietoris topology or Hausdorff topology, instead of Fell topology, then Theorem 13 in [2] does not hold. They provided two examples of totally discrete and locally compact Hausdorff topological ∧-semilattices, in which the map $x \to x^\downarrow$ is not continuous from $(X, \tau, \leqslant)$ to $(C^\downarrow(X), \tau_V, \subseteq)$ and to $(C^\downarrow(X), \tau_H, \subseteq)$, respectively.

In this section, we provide Examples 4.1 and 4.2 of locally compact and order-connected Hausdorff topological ∧-semilattices (all conditions in Theorem 13 in [2] are satisfied). In Examples 4.1, we construct a locally compact and order-connected metric (topological) lattice, in which the map $x \to x^{\downarrow}$ is not continuous at every point from $(X, \tau, \leqslant)$ to $(C^{\downarrow}(X), \tau_V, \subseteq)$. Example 4.2 provides a locally compact and order-connected metric (topological) ∧-semilattice, in which the map $x \to x^{\downarrow}$ is not continuous at every (excepting one) point from $(X, \tau, \leqslant)$ to $(C^{\downarrow}(X), \tau_V, \subseteq)$ and it is not continuous at every point from $(X, \tau, \leqslant)$ to $(C^{\downarrow}(X), \tau_H, \subseteq)$.

**Example 4.1.** Define $X$ to be the following open square in $\mathbb{R}^2$:

$$X = \{(u, v) \in \mathbb{R}^2 : 0 < u, v < 1\}.$$

Let $X$ be equipped with topology $\tau$ induced by the standard ordinary Euclidean norm $\|\cdot\|$ and equipped with the coordinate wise partial order $\leqslant$. Then

(i) $(X, \tau, \leqslant)$ is a locally compact and order-connected Hausdorff topological lattice;
(ii) The canonical map $x \to x^{\downarrow}$ topologically order-embeds $(X, \tau, \leqslant)$ in $(C(X), \tau_F, \subseteq)$;
(iii) The map $x \to x^{\downarrow}$ is not continuous at every point from $(X, \tau)$ to $(C^{\downarrow}(X), \tau_V)$.

*Proof.* The proof of (i) is straightforward and it is omitted here. Then, by part (i) and by Theorem 4.6, or Theorem 13 in [2], part (ii) follows immediately.

Proof of (iii). We first prove that the map $x \to x^{\downarrow}$ is not continuous at every point from $(X, \tau)$ to $(C^{\downarrow}(X), \tau_V)$. For an arbitrary given point $(u_0, v_0)$ in $X$, that is, $0 < u_0, v_0 < 1$, let $l(u_0, v_0)$ be the open segment with none included ending points $(0, v_0)$ and $(u_0, 1)$. Then $l(u_0, v_0)$ is a closed subset of $X$ satisfying $l(u_0, v_0) \cap (u_0, v_0)^{\downarrow} = \emptyset$. It follows that $(u_0, v_0)^{\downarrow} \in (X \setminus l(u_0, v_0))^+$; and therefore, $(X \setminus l(u_0, v_0))^+$ is an $\tau_V$-open neighborhood of $(u_0, v_0)^{\downarrow}$ in $(C^{\downarrow}(X), \tau_V)$. On the other hand, for any point $(u_0, v) \in X$ with $v_0 < v < 1$, $l(u_0, v_0) \cap (u_0, v)^{\downarrow} \neq \emptyset$, that is,

$$(u_0, v)^{\downarrow} \notin (X \setminus l(u_0, v_0))^+, \text{ for any } (u_0, v) \in X \text{ with } v_0 < v < 1.$$

Since $(X \setminus l(u_0, v_0))^+$ is an $\tau_V$-open neighborhood of $(u_0, v_0)^{\downarrow}$, it implies that $(u_0, v)^{\downarrow}$ does not convergent to $(u_0, v_0)^{\downarrow}$, as $v$ goes to $v_0$ (right hand side) with respect to the Vietoris topology $\tau_V$ on the hyperspace $(C^{\downarrow}(X), \tau_V)$. □

The following example is a locally compact and order-connected Hausdorff topological ∧-semilattice, in which the canonical map $x \to x^{\downarrow}$ topologically order-embeds $(X, \tau, \leqslant)$ on to $(C^{\downarrow}(X), \tau_F)$. Meanwhile, this map is not continuous at every point (excepting one) in $X$ with respect to the Vietoris topology $\tau_V$ and not continuous at every point in $X$ with respect to the Hausdorff topology $\tau_H$ on $C^{\downarrow}(X)$.

**Example 4.2.** Define $X$ (a subset (a solid) of $\mathbb{R}^3$) by

$$X = \{(u, v, w) \in \mathbb{R}^3 : -\infty < u, v, w \leq 0 \text{ and } uv + w - 1 \leq 0\}.$$

Let $X$ be equipped with topology $\tau$ induced by the standard ordinary Euclidean norm $\|\cdot\|$ and equipped with the coordinate wise partial order $\leqslant$ on $\mathbb{R}^3$. Then, $(X, \tau, \leqslant)$ is a locally compact and order-connected Hausdorff topological $\wedge$-semilattice. Moreover, the canonical map $x \to x^\downarrow$ satisfies that

(i) It topologically order-embeds $(X, \tau, \leqslant)$ in $(C(X), \tau_F)$;
(ii) It is not continuous at every point from $(X, \tau, \leqslant)$ to $(C^\downarrow(X), \tau_H)$;
(iii) It is not continuous at every point excepting $(0, 0, 0)$ from $(X, \tau, \leqslant)$ to $(C^\downarrow(X), \tau_V)$.

*Proof.* For easier seeing the set $X$, we draw its graph below.

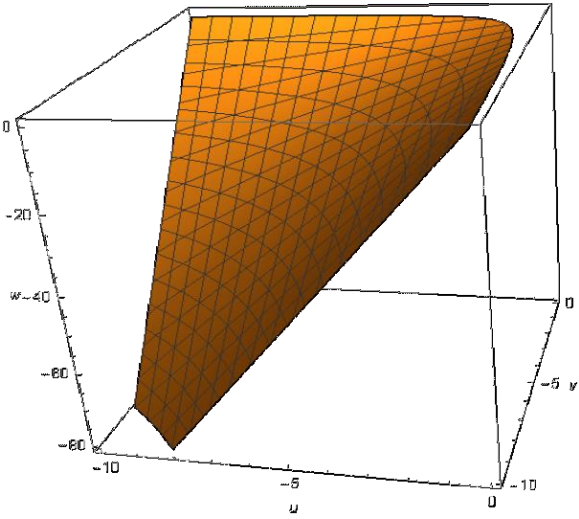

In the setting of the space $X$ as showed in the graph, we see that $X$ is contained in Octant VII (all coordinates are non-positive) in $\mathbb{R}^3$. From geometric view, this solid $X$ of $\mathbb{R}^3$ is enclosed (and included) by the following four surfaces, which are called the faces of $X$.

$S_1 = \{(u, v, w) \in \mathbb{R}^3 : uv + w - 1 = 0, \ -\infty < u, v < 0, \text{ and } -\infty < w \leq 0\}$;
$S_2 = \{(u, 0, w) \in \mathbb{R}^3 : -\infty < u, w \leq 0\}$;
$S_3 = \{(0, v, w) \in \mathbb{R}^3 : -\infty < v, w \leq 0\}$;
$S_4 = \{(u, v, 0) \in \mathbb{R}^3 : -\infty < u, v \leq 0 \text{ and } uv - 1 \leq 0\}$.

$S_1$, $S_2$, $S_3$, and $S_4$ are respectively called the front face, the back face, the right hand side face, and the top face of $X$. More precisely, $S_2$ is quadrant III of the coordinate plane $uOw$, $S_3$ is quadrant III of the coordinate plane $vOw$, and $S_4$ is a closed and unbounded region in the coordinate plane $uOv$. The negative $v$-axis is the right top edge of $X$.

For our convenience, we use the following notations in this proof:

(a) ∧ and ∨ are used in both cases of $(X, \tau, \preccurlyeq)$ and $(\mathbb{R}, \leq)$. For example, for any $u_1, u_2 \in \mathbb{R}$, $u_1 \wedge u_2 = \min\{u_1, u_2\}$.

(b) For every $x = (u, v, w) \in X$, $x^{\downarrow}$ and $x^{\uparrow}$ respectively denote the $\preccurlyeq$-principal ideal and filter in $(X, \tau, \preccurlyeq)$. Since, with the coordinate wise partial order $\preccurlyeq$, $(\mathbb{R}^3, \tau, \preccurlyeq)$ is a locally compact and order-connected topological lattice. Then, for any $x = (u, v, w) \in \mathbb{R}^3$, to distinct the difference between the principal ideals and filters in $(X, \tau, \preccurlyeq)$ and $(\mathbb{R}^3, \tau, \preccurlyeq)$, we use $x^{\downarrow}_{\mathbb{R}^3}$ and $x^{\uparrow}_{\mathbb{R}^3}$ to respectively denote the $\preccurlyeq$-principal ideal and filter in $(\mathbb{R}^3, \tau, \preccurlyeq)$.

At first, we prove that $(X, \tau, \preccurlyeq)$ is a locally compact and order-connected topological ∧-semilattice.

Notice that, as partial of the boundary of $X$, surface $S_1$ is a very smooth surface and $S_2$, $S_3$, and $S_4$ all are partials of planes in $\mathbb{R}^3$; and $\text{int}X = X \setminus (S_1 \cup S_2 \cup S_3 \cup S_4)$ is a connected open subset of $\mathbb{R}^3$. It implies that $(X, \tau)$ is a locally compact subset in $\mathbb{R}^3$. To prove that $(X, \tau, \preccurlyeq)$ is order-connected, take arbitrary two distinct points $x = (u_1, v_1, w_1)$ and $y = (u_2, v_2, w_2)$ in $X$ with $x \preccurlyeq y$, we have

$$x^{\uparrow} \cap y^{\downarrow} = (x^{\uparrow}_{\mathbb{R}^3} \cap y^{\downarrow}_{\mathbb{R}^3}) \cap X.$$

Since $x^{\uparrow}_{\mathbb{R}^3} \cap y^{\downarrow}_{\mathbb{R}^3}$ is either a rectangular box (if $u_1 < u_2$, $v_1 < v_2$ and $w_1 < w_2$), or a rectangle surface (if exact one of the equations hold: $u_1 = u_2$, $v_1 = v_2$ and $w_1 = w_2$), or a closed segment (if exact two of the equations hold: $u_1 = u_2$, $v_1 = v_2$ and $w_1 = w_2$) in $\mathbb{R}^3$, it follows that the subset $x^{\uparrow} \cap y^{\downarrow}$ is the intersection of a rectangular box (may be a rectangle or a segment) and the solid $X$ in $\mathbb{R}^3$. It is indeed a connected subset in $\mathbb{R}^3$.

Next we prove that $(X, \tau, \preccurlyeq)$ is a topological ∧-semilattice. To this ends, arbitrarily take two points $x = (u_1, v_1, w_1)$ and $y = (u_2, v_2, w_2)$ in $X$. We study $x \wedge y$ by the following two cases:

Case 1. If $(u_1 \wedge u_2)(v_1 \wedge v_2) + w_1 \wedge w_2 - 1 \leq 0$, then $(u_1 \wedge u_2, v_1 \wedge v_2, w_1 \wedge w_2) \in X$ and

$$x \wedge y = (u_1 \wedge u_2, v_1 \wedge v_2, w_1 \wedge w_2). \tag{4.1}$$

Case 2. If $(u_1 \wedge u_2)(v_1 \wedge v_2) + w_1 \wedge w_2 - 1 > 0$, then $(u_1 \wedge u_2, v_1 \wedge v_2, w_1 \wedge w_2) \notin X$. In this case, we fixe $u_1 \wedge u_2$ and $v_1 \wedge v_2$ and take $\overline{w} < 0$ such that $(u_1 \wedge u_2)(v_1 \wedge v_2) + \overline{w} - 1 = 0$. From $(u_1 \wedge u_2)(v_1 \wedge v_2) + w_1 \wedge w_2 - 1 > 0$, it follows that

$$\overline{w} = -(u_1 \wedge u_2)(v_1 \wedge v_2) + 1 < w_1 \wedge w_2. \tag{4.2}$$

It implies that $(u_1 \wedge u_2, v_1 \wedge v_2, \overline{w}) \in S_1 \subseteq X$. Since $(u_1 \wedge u_2, v_1 \wedge v_2, w_1 \wedge w_2) \notin X$, considering the space $(\mathbb{R}^3, \tau, \preccurlyeq)$, it follows that

$$(u_1 \wedge u_2, v_1 \wedge v_2, \overline{w}) \preccurlyeq (u_1 \wedge u_2, v_1 \wedge v_2, w_1 \wedge w_2) \text{ in } \mathbb{R}^3.$$

By $(u_1 \wedge u_2, v_1 \wedge v_2, \overline{w}) \in S_1 \subseteq X$, it implies that

$$(u_1 \wedge u_2, v_1 \wedge v_2, \overline{w}) \preccurlyeq x \text{ and } (u_1 \wedge u_2, v_1 \wedge v_2, \overline{w}) \preccurlyeq y \text{ in } X. \tag{4.3}$$

Then we analytically prove that $(u_1 \wedge u_2, v_1 \wedge v_2, \overline{w}) = x \wedge y$. To this ends, for any $(u, v, w) \in X$ with $(u, v, w) \preccurlyeq x$ and $(u, v, w) \preccurlyeq y$, since $\preccurlyeq$ is the coordinate wise order, we must have

$$u \leq u_1 \wedge u_2 \text{ and } v \leq v_1 \wedge v_2.$$

For the fixed $u_1 \wedge u_2$ and $v_1 \wedge v_2$, we know that $(u_1 \wedge u_2, v_1 \wedge v_2, \overline{w}) \in S_1$ and satisfy (4.3). From $(u_1 \wedge u_2)(v_1 \wedge v_2) + \overline{w} - 1 = 0$, it follows that if $w > \overline{w}$, then $(u_1 \wedge u_2)(v_1 \wedge v_2) + w - 1 > 0$. It implies that if $w > \overline{w}$, then $(u_1 \wedge u_2, v_1 \wedge v_2, w) \notin X$. Hence

$$\overline{w} = \max\{w \in \mathbb{R}: (u_1 \wedge u_2, v_1 \wedge v_2, w) \in X\}$$

Hence, for any $(u, v, w) \in X$ with $(u, v, w) \preccurlyeq x$ and $(u, v, w) \preccurlyeq y$, we must have

$$u \leq u_1 \wedge u_2, v \leq v_1 \wedge v_2, \text{ and } w \leq \overline{w}.$$

From $(u_1 \wedge u_2, v_1 \wedge v_2, \overline{w}) \in S_1 \subseteq X$ and by (4.3), the above inequalities imply that

$$x \wedge y = (u_1 \wedge u_2, v_1 \wedge v_2, \overline{w}).$$

Hence, if $(u_1 \wedge u_2)(v_1 \wedge v_2) + w_1 \wedge w_2 - 1 > 0$, then

$$x \wedge y = (u_1 \wedge u_2, v_1 \wedge v_2, \overline{w}) = (u_1 \wedge u_2, v_1 \wedge v_2, -(u_1 \wedge u_2)(v_1 \wedge v_2) + 1). \tag{4.4}$$

(4.1) and (4.4) together provide the expression of $x \wedge y$. Since $\wedge$ is continuous on $\mathbb{R}$ and the function $-(u_1 \wedge u_2)(v_1 \wedge v_2) + 1$ is also continuous, it follows that $\wedge$ is continuous on $X$; and therefore, we proved that $(X, \tau, \preccurlyeq)$ is a topological $\wedge$-semilattice.

Here, for the proof of case 2 that $(X, \tau, \preccurlyeq)$ is a $\wedge$-semilattice, in addition to the above analytical proof, we geometrically show $(u_1 \wedge u_2, v_1 \wedge v_2, \overline{w}) = x \wedge y$ in case 2 such that one can better understand the set $X$ and the proof of case 2.

Let $l$ denote the vertical half line in $\mathbb{R}^3$ with top ending point $(u_1 \wedge u_2, v_1 \wedge v_2, w_1 \wedge w_2)$. $l$ can be considered as the following subset in $\mathbb{R}^3$

$$l = \{(u_1 \wedge u_2, v_1 \wedge v_2, w) \in \mathbb{R}^3: -\infty < w \leq w_1 \wedge w_2\}.$$

Then $l \cap S_1 \neq \emptyset$ and it is a singleton contained in $S_1$. Let $\{(u, v, w)\} = l \cap S_1$. From $(u, v, w) \in l$, it follows that $u = u_1 \wedge u_2$ and $v = v_1 \wedge v_2$. By $(u_1 \wedge u_2, v_1 \wedge v_2, w) \in S_1$, we have $(u_1 \wedge u_2)(v_1 \wedge v_2) + w - 1 = 0$. It follows $w = \overline{w}$, which is exactly the same $\overline{w}$ defined in the proof of case 2 in (4.2). That is

$$l \cap S_1 = \{(u_1 \wedge u_2, v_1 \wedge v_2, \overline{w})\}. \tag{4.5}$$

By the hypothesis in case 2 that $(u_1 \wedge u_2, v_1 \wedge v_2, w_1 \wedge w_2) \notin G$, and noticing that $l$ is the back right edge of $(u_1 \wedge u_2, v_1 \wedge v_2, w_1 \wedge w_2)^{\downarrow}_{\mathbb{R}^3}$ in $\mathbb{R}^3$, from (4.5), we have

$$\{(u, v, w) \in X: (u, v, w) \leqslant (u_1 \wedge u_2, v_1 \wedge v_2, w_1 \wedge w_2)\}$$
$$= (u_1 \wedge u_2, v_1 \wedge v_2, w_1 \wedge w_2)^{\downarrow}_{\mathbb{R}^3} \cap X$$
$$= (u_1 \wedge u_2, v_1 \wedge v_2, \overline{w})^{\downarrow}.$$

It follows immediately that

$$(u_1 \wedge u_2, v_1 \wedge v_2, \overline{w}) = \vee\{(u, v, w) \in X: (u, v, w) \leqslant (u_1 \wedge u_2, v_1 \wedge v_2, w_1 \wedge w_2)\}. \tag{4.6}$$

From $(u_1 \wedge u_2, v_1 \wedge v_2, \overline{w}) \in S_1 \subseteq X$ and by (4.6), it implies

$$x \wedge y = (u_1 \wedge u_2, v_1 \wedge v_2, \overline{w}) = (u_1 \wedge u_2, v_1 \wedge v_2, -(u_1 \wedge u_2)(v_1 \wedge v_2) + 1).$$

Proof of (ii). For an arbitrary fixed point $x_0 \equiv (u_0, v_0, w_0) \in X$, we prove that the canonical map $x \to x^{\downarrow}$ from $(X, \tau, \leqslant)$ to $(C(X), \tau_H)$ is not continuous at $x_0$.

For $x_0 \equiv (u_0, v_0, w_0) \in X$, it holds $u_0 v_0 + w_0 - 1 \leq 0$. Note that the interior of $X$ is an open connected subset of $\mathbb{R}^3$. For the given point $(u_0, v_0, w_0) \in X$, define the following subset of $X$

$$S(w_0) = \{(u, v, w_0) \in X: -\infty < u, v \leq 0 \text{ and } uv + w_0 - 1 \leq 0\}.$$

Then $S(w_0)$ is the level surface of $X$ contained in the level plane $w = w_0$. Since $S(w_0)$ is enclosed by the curves $uv + w_0 - 1 = 0$, negative $u$-exits and negative $v$-exits, it is a polyline connected region (it is not convex). Then, for the fixed $u_0 \leq 0$, $v_0 \leq 0$ and $w_0 \leq 0$, there is a point $x_1 \equiv (u_1, v_1, w_0) \in X$ satisfying $u_1 v_1 + w_0 - 1 \leq 0$, for some $u_1 < u_0 \leq 0$ and $0 \geq v_1 \geq v_0$, such that the whole closed segment with ending points $x_1$ and $x_0$ is contained in $S(w_0) \subseteq X$. That is, $S(w_0)$ contains all the linear combinations $x_\alpha$ of $x_0$ and $x_1$, where $x_\alpha$ is defined as

$$x_\alpha = \alpha x_1 + (1-\alpha)x_0 = (\alpha u_1 + (1-\alpha) u_0, \alpha v_1 + (1-\alpha) v_0, w_0) \equiv (u_\alpha, v_\alpha, w_0), \text{ for } \alpha \in (0, 1),$$

By $u_1 < u_0 \leq 0$, $0 \geq v_1 \geq v_0$ and $u_\alpha \equiv \alpha u_1 + (1-\alpha)u_0$, $v_\alpha \equiv \alpha v_1 + (1-\alpha)v_0$, we have

$$u_1 < u_\alpha < u_0 \leq 0 \text{ and } v_0 \leq v_\alpha \leq v_1 \leq 0, \text{ for any } \alpha \in (0, 1). \tag{4.7}$$

Then, $u_\alpha v_\alpha + w_0 - 1 \leq 0$, that is,

$$x_\alpha \equiv (u_\alpha, v_\alpha, w_0) \in S(w_0) \subseteq X, \text{ for any } \alpha \in (0, 1).$$

In the setting of the space $X$ and the coordinates of $\mathbb{R}^3$ as showed in the graph, $x_0^{\downarrow}$ has back top edge and right top edge. We randomly choose the right top edge of $x_0^{\downarrow}$ to use in this proof (that is, it can be similarly proved if we choose the back top edge of $x_0^{\downarrow}$).

Let $l(x_0) \equiv l(u_0, v_0, w_0)$ be the right top edge of $x_0^{\downarrow}$. For $x_0 \equiv (u_0, v_0, w_0)$ with $u_0 \leq 0$, we study $l(x_0)$ by two cases: $u_0 = 0$ and $u_0 < 0$.

If $u_0 = 0$, then, $x_0 \equiv (0, v_0, w_0) \in S_3 \subseteq X$. In this case, the right top edge of $x_0^\downarrow$ is a half line in $S_3$ with ending point $x_0 \equiv (0, v_0, w_0)$ ($S_3$ is quadrant III of the coordinate plane $vOw$) and $x_0^\downarrow$ is considered as the following subset of $S_3 \subseteq X$

$$l(x_0) = \{(0, v, w_0) \in S_3 : v_0 \geq v > -\infty\}. \tag{4.8}$$

If $u_0 < 0$, then $x_0^\downarrow \cap S_3 = \emptyset$. In this case, the (quarter) plane $u = u_0$ named by $P(u_0)$ is considered as the following subset of $\mathbb{R}^3$

$$P(u_0) = \{(u_0, v, w) \in \mathbb{R}^3 : -\infty < v, w \leq 0\}.$$

Since $u_0 < 0$, $P(u_0)$ and $S_3$ are parallel and disjoint. In this case, as the right top edge of $x_0^\downarrow$, $l(x_0)$ is a polyline in $P(u_0)$. In order to find the equation of $l(x_0)$, solving for $v$ from $u_0 v + w_0 - 1 = 0$, we obtain $v_2 = \frac{1 - w_0}{u_0}$. It satisfies $u_0 v_2 + w_0 - 1 = 0$; and therefore, $(u_0, v_2, w_0) \in S_1$. By $u_0 v_0 + w_0 - 1 \leq 0$, $u_0 v_2 + w_0 - 1 = 0$, and $u_0 < 0$, it implies $v_0 \geq v_2$. Hence, if $u_0 < 0$, as the right top edge of $x_0^\downarrow$, this polyline $l(x_0)$ is contained in $P(u_0)$ and is considered as the following subset of $X$

$$l(x_0) = \begin{cases} \{(u_0, v, w_0) \in X : v_0 \geq v \geq v_2\}, \\ \{(u_0, v, 1 - u_0 v) \in X : v_2 \geq v > -\infty\}. \end{cases} \tag{4.9}$$

(4.8) and (4.9) together provide the expressions of $l(x_0)$, which is the right top edge of $x_0^\downarrow$. We see that, from the solution $v_2$ of (4.9), when $u_0 \to 0^-$, $v_2 = \frac{1 - w_0}{u_0} \to -\infty$; and therefore (4.9) turns to be (4.8). Hence, (4.8) can be considered as the special case of $u_0 = 0$.

For any given $\alpha \in (0, 1)$, $x_\alpha \equiv (u_\alpha, v_\alpha, w_0) \in X$, that is, $u_\alpha v_\alpha + w_0 - 1 \leq 0$. Let $l(x_\alpha) \equiv l(u_\alpha, v_\alpha, w_0)$ be the right top edge of $x_\alpha^\downarrow$. From (4.7), for the given $u_\alpha$, $v_\alpha$ satisfying $u_1 < u_\alpha < u_0 \leq 0$, $v_1 \geq v_\alpha \geq v_0$, and already fixed $w_0 \leq 0$, let $P(u_\alpha)$ name the (quarter) plane $u = u_\alpha$, which is considered as the following subset of $\mathbb{R}^3$

$$P(u_\alpha) = \{(u_\alpha, v, w) \in \mathbb{R}^3 : -\infty < v, w \leq 0\}.$$

$P(u_\alpha)$, $P(u_0)$ and $S_3$ are parallel and, by $u_\alpha < u_0 \leq 0$, we have

$$P(u_\alpha) \cap P(u_0) = \emptyset, \text{ and } P(u_\alpha) \cap S_3 = \emptyset, \text{ for } u_1 < u_\alpha < u_0 \leq 0.$$

The right top edge $l(x_\alpha)$ of $x_\alpha^\downarrow$ is a polyline in the (quarter) plane $P(u_\alpha)$. To find the equation of $l(x_\alpha)$, solving for $v$ from $u_\alpha v + w_0 - 1 = 0$, we obtain $v_3 = \frac{1 - w_0}{u_\alpha}$. It satisfies $u_\alpha v_3 + w_0 - 1 = 0$; and therefore, $(u_\alpha, v_3, w_0) \in S_1$. By $u_\alpha v_\alpha + w_0 - 1 \leq 0$, $u_\alpha v_3 + w_0 - 1 = 0$, and $u_\alpha < 0$, it implies $v_\alpha \geq v_3$. Then, as the right top edge of $x_\alpha^\downarrow$, $l(x_\alpha)$ is a polyline in the plane $P(u_\alpha)$ considered as the following subset of $X$

$$l(x_\alpha) = \begin{cases} \{(u_\alpha, v, w_0) \in X: v_\alpha \geq v \geq v_3\} \\ \{(u_\alpha, v, 1 - u_\alpha v) \in X: v_3 \geq v > -\infty\} \end{cases}. \tag{4.10}$$

If $u_0 = 0$, then $l(x_0)$ is a half line in $S_3$ that has slope 0. In this case, for any given $\alpha \in (0, 1)$, $l(x_\alpha)$ is a polyline in the (quarter) plane $P(u_\alpha)$ with a positive slope $- u_\alpha$ on $(-\infty, v_3]$. Hence $l(x_\alpha)$ is on the left and below of $l(x_0)$. If $u_\alpha < u_0 \leq 0$, from $v_3 = \frac{1-w_0}{u_\alpha}$ and $v_2 = \frac{1-w_0}{u_0}$, it follows that $v_3 > v_2$. In this case, on $(-\infty, v_2]$, $l(x_0)$ is a polyline in the (quarter) plane $P(u_0)$ with slope $-u_0 > 0$, and $l(x_\alpha)$ is a polyline in the (quarter) plane $P(u_\alpha)$ with a positive slope $- u_\alpha$ on $(-\infty, v_3]$. Since $-u_\alpha > -u_0$, it follows that $l(x_\alpha)$ is on the left and below of $l(x_0)$. It concludes that, for any given $\alpha \in (0, 1)$, the right top edge $l(x_\alpha)$ of $x_\alpha^\downarrow$ (a polyline in the (quarter) plane $P(u_\alpha)$) is on the left and below of the right top edge $l(x_0)$ of $x_0^\downarrow$. It implies that, for any $z \equiv (u_0, v, w) \in l(x_0)$,

$$d(z, x_\alpha^\downarrow) = d(z, l(x_\alpha)).$$

Let $l_P(x_\alpha, u_0)$ be the projection of $l(x_\alpha)$ to the (quarter) plane $P(u_0)$, which is parallel to both $P(u_\alpha)$ and $S_3$ satisfying $P(u_0) \cap P(u_\alpha) = \emptyset$, for $u_1 < u_\alpha < u_0 \leq 0$. Then, from (4.10), $l_P(x_\alpha, u_0)$ is also a polyline in the (quarter) plane $P(u_0)$ satisfying

$$l_P(x_\alpha, u_0) = \begin{cases} \{(u_0, v, w_0) \in X: v_\alpha \geq v \geq v_3\} \\ \{(u_0, v, 1 - u_\alpha v) \in X: v_3 \geq v > -\infty\} \end{cases},$$

where, $\alpha \in (0, 1)$, $u_1 < u_\alpha < u_0 \leq 0$ and $0 \geq v_1 \geq v_\alpha \geq v_0$.

From geometric view, in the same (quarter) plane $P(u_0)$, when $v_3 > v_2 > v > -\infty$, the two half lines $l_P(x_\alpha, u_0)$ and $l(x_0)$ have slopes $-u_\alpha$ and $u_0$, respectively. Since $-u_\alpha > -u_0$, we must have

$$H(l(x_\alpha), l(x_0)) \geq H(l_P(x_\alpha, u_0), l(x_0)) = \infty. \tag{4.11}$$

Here, we sketch the details of the proof of (4.11). Since $u_\alpha < u_0 \leq 0$, and $l(x_\alpha)$ is on the left and below of $l(x_0)$, it follows that, $l_P(x_\alpha, u_0)$ is below $l(x_0)$ in the same plane $P(u_0)$. Take an arbitrary point $z \equiv (u_0, v', w') \in l(x_0)$, for some $v_0 \geq v' > -\infty$, we have

$$z = (u_0, v', w') = \begin{cases} (0, v', w_0), & \text{if } u_0 = 0, \\ (u_0, v', w_0), & \text{if } u_0 < 0 \text{ and } v_3 \leq v' \leq v_0, \\ (u_0, v', 1 - u_0 v'), & \text{if } u_0 < 0 \text{ and } v' \leq v_3. \end{cases} \tag{4.12}$$

Since $P(u_\alpha) \cap P(u_0) = \emptyset$, and $P(u_\alpha) \cap S_3 = \emptyset$, for $u_1 < u_\alpha < u_0 \leq 0$, and $l_P(x_\alpha, u_0)$ is the projection of $l(x_\alpha)$ from $P(u_\alpha)$ to $P(u_0)$, it implies that, for $z = (u_0, v', w') \in l(x_0) \subseteq P(u_0)$, we have

$$d(z, l(x_\alpha)) \geq d(z, l_P(x_\alpha, u_0)).$$

$l(x_\alpha)$ is the right top edge of $x_\alpha^\downarrow$ and it is on the left below of $l(x_0)$, it follows that

$$d(z, x_\alpha^\downarrow) = d(z, l(x_\alpha)) \geq d(z, l_P(x_\alpha, u_0)).$$

Now we calculate $d(z, l_P(x_\alpha, u_0))$. In the case that $u_\alpha < u_0 = 0$, $l_P(x_\alpha, u_0)$ and $l(x_0)$ coincide in the part $v_3 \leq v' \leq v_\alpha \leq v_0$. It implies that $d(z, l_P(x_\alpha, u_0)) = 0$, for $z = (u_0, v', w_0) \in l(x_0)$ with $v_3 \leq v' \leq v_\alpha \leq v_0$. In the case that $u_\alpha < u_0 \leq 0$, from $v_3 = \frac{1-w_0}{u_\alpha}$ and $v_2 = \frac{1-w_0}{u_0}$, it follows that $v_3 > v_2$. It implies that $l_P(x_\alpha, u_0)$ and $l(x_0)$ also coincide in the part $v_3 \leq v' \leq v_\alpha \leq v_0$. It implies that $d(z, l_P(x_\alpha, u_0)) = 0$, for $z = (u_0, v', w_0)$ in the case $u_\alpha < u_0 \leq 0$ and $v_3 \leq v' \leq v_\alpha \leq v_0$.

Since we will consider the situation that $v' \to -\infty$, we may only consider the case $v' < v_2 < v_3$. Since $u_\alpha < 0$, when $v < v_2 (< v_3)$, the line $l_P(x_\alpha, u_0)$ is in the plane $P(u_0)$ and has equation $w = -u_\alpha v + 1$. For the part $v_3 > v_2 > v' > -\infty$, to calculate $d(z, l_P(x_\alpha, u_0))$, we draw a line in the (quarter) plane $P(u_0)$ that passes through point $z = (u_0, v', w') \in l(x_0)$ and is perpendicular to the line $w = -u_\alpha v + 1$. Then this line has equation $w = \frac{1}{u_\alpha}(v - v') + w'$. Then these two lines $w = \frac{1}{u_\alpha}(v - v') + w'$ and $w = -u_\alpha v + 1$ have the joint point in $P(u_0)$

$$z_\alpha \equiv \left(u_0, \frac{v' + u_\alpha - u_\alpha w'}{1 + u_\alpha^2}, 1 - \frac{u_\alpha(v' + u_\alpha - u_\alpha w')}{1 + u_\alpha^2}\right).$$

Calculate $d(z, z_\alpha)$, we obtain

$$d(z, z_\alpha) = \frac{|1 - u_\alpha v' - w'|}{\sqrt{1 + u_\alpha^2}}, \text{ if } v_3 > v_2 > v' > -\infty. \tag{4.13}$$

From (4.12), if $v_3 > v_2 > v' > -\infty$, we have

$$z = (u_0, v', w') = \begin{cases} (0, v', w_0), & \text{if } u_0 = 0, \\ (u_0, v', 1 - u_0 v'), & \text{if } u_0 < 0 \text{ and } v' < v_2 \end{cases}$$

Substituting it into (4.13), we obtain

$$d(z, z_\alpha) = \begin{cases} \frac{|1 - u_\alpha v' - w'|}{\sqrt{1 + u_\alpha^2}}, & \text{if } u_0 = 0, \\ \frac{|v'||u_\alpha - u_0|}{\sqrt{1 + u_\alpha^2}}, & \text{if } u_0 < 0 \text{ and } v' < v_3. \end{cases} \tag{4.14}$$

Since $\frac{|u_\alpha - u_0|}{\sqrt{1 + u_\alpha^2}} > 0$, and $u_\alpha < 0$, from (4.14), it follows that, for any given $x_\alpha \equiv (u_\alpha, v_\alpha, w_0) \in X$ with $\alpha \in (0, 1)$, we have

$$d(z, l_P(x_\alpha, u_0)) \to \infty, \text{ as } v' \to -\infty. \tag{4.15}$$

By (4.15), it implies that

$$H(x_\alpha^\downarrow, x_0^\downarrow) \geq \sup\{d(z, l(x_\alpha)): z \in l(x_0)\} \geq \sup\{d(z, l_P(x_\alpha, u_0)): z \in l(x_0)\} = \infty. \tag{4.16}$$

Then, for this arbitrary fixed $x_0 \equiv (u_0, v_0, w_0) \in X$, by (4.16), we have

$$H(x_\alpha^\downarrow, x_0^\downarrow) = \infty, \text{ for } x_\alpha \equiv (u_\alpha, v_\alpha, w_0) \in S(w_0) \subseteq X \text{ with } \alpha \in (0, 1).$$

Since $(u_\alpha, v_\alpha, w_0) \to (u_0, v_0, w_0)$, as $\alpha \to 0$, it implies that the canonical map $x \to x^\downarrow$ from $(X, \tau)$ to $(C(X), \tau_H)$ is not continuous at this arbitrary given point $x_0 \equiv (u_0, v_0, w_0) \in X$. It proves part (ii).

Proof of (iii). We prove that the canonical map $x \to x^\downarrow$ is not continuous at every point in $X \setminus \{(0, 0, 0)\}$ from $(X, \tau, \preccurlyeq)$ to $(C(X), \tau_V)$.

Take an arbitrary point $x_0 \equiv (u_0, v_0, w_0) \in X$ with $(u_0, v_0, w_0) \neq (0, 0, 0)$.

Case 1. Either $u_0 < 0$, or $v_0 < 0$, or both and for any $w_0 \leq 0$. Suppose that $v_0 < 0$. Then, we have

$$(u_0, v, w) \in X, \text{ for any } v \in [v_0, 0] \text{ and } -\infty < w \leq w_0.$$

For $(u_0, v_0, w_0) \in X$ with $v_0 < 0$, let $S[v_0, 0]$ denote the following subset of $X$

$$S[v_0, 0] \equiv \{(u_0, v, w): v \in [v_0, 0] \text{ and } w \leq w_0\} \subseteq X.$$

Then $S[v_0, 0]$ is a rectangular section and it is unbounded from below. Define a curve $C(v_0, 0]$ as follows:

$$C(v_0, 0] \equiv \{(u_0, v, \frac{v}{v-v_0} + w_0): v \in (v_0, 0]\}.$$

From $\frac{v}{v-v_0} \leq 0$, for $v \in (v_0, 0]$, it follows that $C(v_0, 0] \subseteq S[v_0, 0] \subseteq X$. Denote the vertical half-line with ending point $(u_0, v, w_0)$ by

$$L(v) \equiv \{(u_0, v, w): w \leq w_0\}, \text{ for } v \in (v_0, 0].$$

The vertical half-lines $L(v)$ have that following properties:

(a) $L(v) \subseteq S[v_0, 0]$, for every $v \in (v_0, 0]$;
(b) $L(v) \cap C(v_0, 0] \neq \emptyset$, and it is a singleton, for every $v \in (v_0, 0]$;
(c) $L(v_0) \cap C(v_0, 0] = \emptyset$.

$C(v_0, 0]$ is a closed subset in $X$. Since $L(v_0)$ is the vertical right back edge of $(u_0, v_0, w_0)^\downarrow$, from the above property (c), it implies that $x_0^\downarrow \cap C(v_0, 0] = \emptyset$; and therefore, $(X \backslash C(v_0, 0])_{\tau_V}^+$ is an $\tau_V$-open neighborhood of $x_0^\downarrow$ in $(C(X), \tau_V)$. On the other hand, for every fixed $v \in (v_0, 0]$, $L(v)$ is the vertical right back edge of $(u_0, v, w_0)^\downarrow$. From above property (b), it implies that,

$$(u_0, v, w_0)^\downarrow \cap C(v_0, 0] \supseteq L(v) \cap C(v_0, 0] \neq \emptyset, \text{ for every } v \in (v_0, 0].$$

That is,

$$(u_0, v, w_0)^\downarrow \notin (X \backslash C(v_0, 0])_{\tau_V}^+, \text{ for every } v \in (v_0, 0].$$

It implies that the canonical map $x \to x^\downarrow$ from $(X, \tau)$ to $(C(X), \tau_V)$ is not continuous at this arbitrary given point $x_0 \equiv (u_0, v_0, w_0)$ with $v_0 < 0$. We can similarly prove that map $x \to x^\downarrow$ is not $\tau - \tau_V$ continuous at $x_0 \equiv (u_0, v_0, w_0)$ with $u_0 < 0$.

Case 2. $u_0 = 0$, $v_0 = 0$ and $w_0 < 0$. Then, from $(0, 0, w) \in X$, for any $w \in [w_0, 0]$, we similarly define a subset $s[w_0, 0]$ of $X$ by

$$s[w_0, 0] \equiv \{(0, v, w): v \in (-\infty, 0] \text{ and } w \in [w_0, 0]\} \subseteq X.$$

$s[w_0, 0]$ is a rectangular section unbounded from right. Define a curve $c(w_0, 0]$ in $s[w_0, 0]$ as follows:

$$c(w_0, 0] \equiv \{(0, \frac{w}{w - w_0}, w): \text{ for } w \in (w_0, 0]\}.$$

From $\frac{w}{w - w_0} \leq 0$, for every $w \in (w_0, 0]$, it follows that $c(w_0, 0] \subseteq s[w_0, 0]$. For any given fixed $w \in [w_0, 0]$, denote a half-line with ending point $(0, 0, w)$ by

$$l(w) \equiv \{(0, v, w): v \in (-\infty, 0]\}.$$

These half-lines $l(w)$, for $w \in [w_0, 0]$, have that following properties:

(d) $l(w) \cap c(w_0, 0] \neq \emptyset$ and it is a singleton, for every $w \in (w_0, 0]$;
(e) $l(w_0) \cap c(w_0, 0] = \emptyset$.

$c(w_0, 0]$ is a closed subset in $X$. Since $l(w_0)$ is the top right edge of $(0,0, w_0)^\downarrow$, from the above property (e), it implies that $x_0^\downarrow \cap c(w_0, 0] = \emptyset$; and therefore, $(X \backslash c(w_0, 0])_{\tau_V}^+$ is an $\tau_V$-open neighborhood of $x_0^\downarrow$ in $(C(X), \tau_V)$.

On the other hand, for every $w \in (w_0, 0]$, since $L(w)$ is the top right edge of $(0,0, w)^\downarrow$, from above property (d), it implies

$$(0,0, w_0)^\downarrow \cap c(w_0, 0] \supseteq l(w) \cap c(w_0, 0] \neq \emptyset, \text{ for every } w \in (w_0, 0].$$

That is,
$$(0,0,w)^{\downarrow} \notin (X\backslash c(w_0, 0])^{+}_{\tau_V}, \text{ for every } w \in (w_0, 0].$$

It implies that the canonical map $x \to x^{\downarrow}$ from $(X, \tau)$ to $(C(X), \tau_V)$ is not continuous at this arbitrary given point $x_0 \equiv (0, 0, w_0)$ with $w_0 < 0$.

Finally, we prove that the canonical map $x \to x^{\downarrow}$ is continuous at $(0, 0, 0)$ in $X$ from $(X, \tau)$ to $(C(X), \tau_V)$.

Since $(0,0,0)^{\downarrow} = X$, which contains every closed subset of $X$, it follows that every base $\tau_V$-open neighborhood of $(0,0,0)^{\downarrow}$ must be the type of $O^{-}$, for some nonempty open subset $O \subseteq X$, and vice versa. Hence, for an arbitrary base $\tau_V$-open neighborhood $O^{-}$ of $(0,0,0)^{\downarrow}$, for some nonempty open subset $O \subseteq X$, the openness of $O$ implies that there is a point $z_0 \equiv (u_0, v_0, w_0) \in O$ with $\min\{|u_0|, |v_0|, |w_0|\} > 0$. Take any $0 < r < \min\{|u_0|, |v_0|, |w_0|\}$. Let $B(r) = \{x = (u, v, w) \in X : 0 \geq u, v, w > -r\}$. Then $B(r)$ is an $\tau$-open neighborhood of $(0, 0, 0)$. It satisfies that, for any $x \in B(r)$, $z \in x^{\downarrow} \cap O$, that is $x^{\downarrow} \cap O \neq \emptyset$; and therefore

$$x^{\downarrow} \in O^{-}, \text{ for any } x \in B(r).$$

Hence the canonical map $x \to x^{\downarrow}$ is continuous at $(0, 0, 0)$ from $(X, \tau)$ to $(C(X), \tau_V)$. □

## Acknowledgments

The author is very grateful to Professor Efe A. Ok for introducing me to their valuable article [2]. The key ideas behind the constructions of the examples in this article are based on the techniques developed in the article [2]. The author also sincerely thanks Professor Ok for very valuable communications and for many useful answers of my questions for the accomplishment of this article.